\documentclass[a4paper,11pt]{article} %

\usepackage{amsthm}
\usepackage{amssymb}
\usepackage{graphicx,amsmath} 
\usepackage{tikz}
\usepackage{natbib}
\usepackage{graphics}
\usepackage{subcaption}
\usepackage{setspace}
\usepackage[margin=3cm]{geometry}
\usepackage{titlesec}
\usepackage{microtype}
\usepackage{breakcites}
\usepackage[colorlinks=true,linkcolor=blue, urlcolor=blue, citecolor=blue]{hyperref}

\titleformat{\section}{\filcenter \scshape}{\thesection.}{.5em}{}
\titleformat{\subsection}{\filcenter \small}{\thesubsection}{.5em}{}

\newtheorem{theorem}{Theorem}

\theoremstyle{definition}
\newtheorem{definition}{Definition}

\begin{document}

\title{Playing Tennis without Envy}
\author{Josu\'e Ortega\\University of Glasgow\\ Forthcoming in Mathematics Today}

\date{}
\maketitle

\begin{abstract}
A group of friends organize their tennis games by submitting each their availability over the weekdays. They want to obtain an assignment such that: each game must be a double tennis match, i.e. requires four people, and nobody plays in a day he is unavailable. Can we construct assignments that will always produce efficient, fair, and envy-free outcomes? The answer is no, and extends to any sport that requires any group size.  
\end{abstract}
\thispagestyle{empty}

In the June 2016 edition of the magazine Mathematics Today, \citet{maher2016} described an algorithm to assign tennis double matches among his circle of friends. The algorithm takes as input the players’ availability for the weekdays, and maximizes the number of tennis games, subject to three constraints: 1) no agent plays more than once per day, 2) each match has exactly four players, and 3) no agent plays on a day he is not available.

The algorithm solves a linear program to maximize the number of games achieved but the solution is generally not unique. Hence Maher selects among those the ones that maximize the number of players that get at least one game. If several assignments remain, then he chooses the ones that maximize the number of players that get at least two games. In case uniqueness is yet not achieved, he selects one solution randomly among those, which is the one implemented and communicated to each player. Maher writes ``{\it the (members of the group) appear to trust in the fairness and efficiency of the algorithm}''.

It is clear that the final assignment is efficient as it maximizes the number of matches. But is it really fair? In the preferences that appear in Maher's article, presented in Table 1 below for convenience, George StC declares to be available for 4 days - he is the most flexible player as can play basically any day. However, he only gets one match: the final assignment appears in brackets in Table \ref{tab:tab1}. Six players (Barry T, Peter W, Colin C, Keith B, Brian F, and Peter K) were all half as flexible as George StC and got twice as games as him. George StC may argue the final assignment is treating him unfairly, and probably most readers would agree with him. Furthermore, everybody except Gordon B has an assignment at least as good as George StC.\footnote{This is a actually a simplification of Maher's problem in which we do not consider individual quotas.}

\begin{table}[ht]
\centering
\caption{Player's availability from Prof. Maher's tennis group.}
\label{tab:tab1}
\begin{tabular}{|l|l|l|l|l|l|l|}
\hline
Names      & Mon   & Tues  & Wed   & Thurs & Fri   & Times \\
\hline
Barry T    & 0     & 0     & 1 (1) & 1 (1) & 0     & 2 (2) \\
Tom B      & 1 (0) & 1 (1) & 0     & 1 (1) & 0     & 3 (2) \\
Gordon B   & 0     & 0     & 0     & 0     & 1 (0) & 1 (0) \\
Peter W    & 1 (1) & 1 (1) & 0     & 0     & 0     & 2 (2) \\
Colin C    & 1 (1) & 0     & 0     & 1 (1) & 0     & 2 (2) \\
Mike M     & 0     & 1 (1) & 1 (0) & 1 (1) & 1 (0) & 4 (2) \\
Keith I    & 0     & 1 (1) & 1 (1) & 0     & 0     & 2 (1) \\
Alan C     & 1 (0) & 0     & 0     & 1 (1) & 0     & 2 (1) \\
John S     & 0     & 1 (1) & 0     & 0     & 0     & 1 (1) \\
Keith B    & 1 (1) & 0     & 1 (1) & 0     & 0     & 2 (2) \\
George StC & 1 (0) & 1 (0) & 1 (0) & 1 (1) & 0     & 4 (1) \\
Michael L  & 0     & 0     & 1 (1) & 0     & 0     & 1 (1) \\
Phil M     & 0     & 1 (1) & 0     & 0     & 0     & 1 (1) \\
Brian F    & 1 (1) & 1 (1) & 0     & 0     & 0     & 2 (2) \\
Peter K    & 0     & 1 (1) & 0     & 1 (1) & 0     & 2 (2) \\
Willie McM & 0     & 0     & 0     & 1 (1) & 0     & 1 (1) \\
Ken L      & 0     & 1 (1) & 0     & 0     & 0     & 1 (1) \\
\hline
Total      & 7     & 10    & 6     & 8     & 2     &    \\
\hline
\end{tabular}
\end{table}

The property we described is a variant of game-theoretic envy-freeness (see \citet{moulin1995}). This is, in the assignment presented in Table 1, George StC is envious of Barry T, who was less flexible but got more games, or envious of Peter W, because George StC prefers his assignment. An algorithm that always produces envy-free assignment has an important property: players do not want to fictitiously reduce their availability in order to get more games. Just by one player misreporting his true availability, the assignment described previously could change dramatically.

In this tennis assignment problem, that we describe formally below, players have dichotomous preferences over the days: either they want to play or they do not. These preferences were first studied by \citet{bogomolnaia2004}, and the preferences here represent a natural extension of those: agents want to play in as many feasible days as possible. However, any assignments that gives them a game on a day they are not available is considered worse than having no games at all. Hence, players' preferences can be captured with a subset of all possible days. The constraints that 4 people are required for a game has been previously imposed by \citet{shubik1971} over assignments of one day only. This note is the natural extension of these two environments.

\section{Model}

Let $A^*$ be a $n\times m$ binary matrix containing the preferences of each person $i=(1,\ldots,n)$ about playing on day $j=(1,\ldots,m)$; the entry $a^*_{ik}=1$ if person $i$ is available to play on day $k$ and 0 otherwise. $A^*$ will be called a tennis problem and represents the players' preferences, who are indifferent about their game partners and just care about the days in which they play.

A matrix $A^*$ can be reduced to a matrix $A$ by deleting all days when there are not enough people available to create even one match, as in Friday in Table \ref{tab:tab1}. A further reduction can be performed by eliminating people that are not available on any remaining days. The days and players which are eliminated are irrelevant for the type of solutions we will consider, and hence we will work from now on with the corresponding irreducible tennis problem $A$. Formally, an irreducible problem $A$ satisfies:
\begin{eqnarray}
\forall i \in \{1,\ldots,n\},&\quad& \sum_{k=1}^{m} a_{ik} \geq 1\\
\forall k \in \{1,\ldots,m\},&\quad & \sum_{i=1}^n a_{ik} \geq 4
\end{eqnarray}

A solution to $A$ is a binary matrix $X(A)$, whose elements have the same interpretation as in $A$, satisfying the following constraints:
\begin{eqnarray}
\label{eq:ir}\forall i \in\{1,\ldots,n\}, \forall k \in \{1,\ldots,m\}, && a_{ik}=0 \implies  x_{ik}=0\\
\label{eq:feasibility}\forall k \in \{1,\ldots,m\}, && \left(\sum_i x_{ik}\right)\mod 4=0 
\end{eqnarray}

There are three types of conditions we look for: efficiency, fairness, and strong envy-freeness. We look at them in that order.

\begin{definition}An assignment $X$ is {\bf efficient} if there is no assignment $X'$ such that $\sum_i \sum_k x'_{ik}> \sum_i \sum_k x_{ik}$. \end{definition}

Efficient assignments are exactly those that are Pareto optimal, i.e. those in which no player can be made better off without hurting another one. The next one property considers the games received by the individual in the society who is in the worst position, then the ones received by the second worst one, and so on, in the spirit of John Rawls' leximin criterion.

\begin{definition}Let $G_q(X)$ denote the number of rows in $X$ such that $\sum_k x_{ik}\geq q$, for any integer $q$: this is the number of players with at least $q$ games. An assignment $X$ is fairer than another assignment $X'$ if there exists an integer $q$ for which $G_q(X)>G_q(X')$ and for any integer $q'<q$, $G_q(X)=G_q(X')$. An assignment $X$ is {\bf fair} if there is no other assignment which is fairer. \end{definition}

This notion of fairness implies an optimality condition: while we can construct assignments that are efficient but not fair, every fair assignment is efficient (otherwise another match could be created giving some people more matches, contradicting the fairness property). Finally we have a variant of envy-freeness.

\begin{definition}An assignment $X$ is {\bf strongly envy-free} if for any two players $i,j$ with $\sum_k a_{ik}>\sum_j a_{jk}$, we have $\sum_k x_{ik} \geq \sum_k x_{jk}$. \end{definition}

Strong envy-freeness captures the idea that more flexible people should not be penalized by the assignment. Strong envy cannot arise from days on which only the envious person is available to play, as we are working with the corresponding irreducible tennis problem. It is called strong because standard envy-freeness means that nobody prefers someone's else schedule, a property which is clearly too hard to satisfy in this case.

\section{An Impossibility Result}

There are tennis problems that admit no solutions that is strongly envy-free and efficient.

\begin{table}[ht]
\centering
\caption{The impossibility of strong EF and efficient assignments.}
\label{tab:tab3}
\begin{tabular}{|l|l|l|l|l|l|l|}
\hline
Names   & Mon & Tues & Wed & Thur & Frid & Times \\
\hline
$a$     & 1   & 1    & 0   & 0     & 0   & 2     \\
$b$     & 1   & 1    & 0   & 0     & 0   & 2     \\
$c$		 	& 1   & 1    & 0   & 0     & 0   & 2     \\
$d$		  & 1   & 1    & 0   & 0     & 0   & 2     \\
$e$	    & 0   & 0    & 1   & 1     & 1   & 3     \\
$f$	    & 0   & 0    & 1   & 1     & 1   & 3     \\
$g$ 	  & 0   & 0    & 1   & 1     & 1   & 3     \\
$h$	    & 0   & 0    & 1   & 1     & 1   & 3     \\
$i$	    & 0   & 0    & 1   & 1     & 1   & 3     \\
$j$	    & 0   & 0    & 1   & 1     & 1   & 3     \\
$k$		  & 0   & 0    & 1   & 1     & 1   & 3     \\
\hline
Total   & 4   & 4    & 7   & 7     & 7   &    \\
\hline
\end{tabular}
\end{table}

In the tennis assignment problem in Table \ref{tab:tab3} there are eleven players. Note that by efficiency, we need to organize five games, one each day. This implies that players $a$, $b$, $c$, and $d$ get two games. Then, for whatever way we assign the remaining players to the games on Wednesday, Thursday, and Friday, one of the agents’ $f$ to $k$ gets at most one game while he has an availability of three, so he is envious of any player with availability two. This shows the aforementioned impossibility.

One may think of assignments for other sports. For example, a game of poker that requires three players exactly. Or an indoor football match that requires 10 players. In general, let a $q$-sport assignment problem be the one that requires $q$ agents per day, with the tennis assignment problem being its particular case when $q=4$. A simple modification of the example in Table \ref{tab:tab3} shows that our previous conclusion generalizes for arbitrary $q$-sport assignment problems (although for $q=2$ one needs to add more days). This is

\begin{theorem}
For any integer $q \geq 2$, there exists $q$-sport assignment problems which have no solution that is efficient and strongly envy-free (henceforth no solution that is fair and strongly envy-free).
\end{theorem} 

While we obtained a negative result, we leave many questions unanswered regarding how to construct optimal tennis assignments. We note that this assignment problem, despite being very simple, is close in spirit to the stable marriage problem proposed by \cite{gale1962}, which has led to the improvement of real-life assignments such as those between colleges and students, organs and donors, or junior doctors and hospitals, and for which Shapley received the Nobel Prize in Economics in 2012. Hence, this type of problem, while simple, is always worth considering.

\setlength{\bibsep}{0.1cm}

\end{document}